\documentclass[
aps,%
12pt,%
final,%
notitlepage,%
oneside,%
onecolumn,%
nobibnotes,%
nofootinbib,%
superscriptaddress,%
noshowpacs,%
centertags]%
{revtex4}
\usepackage{amsmath,amssymb}
\usepackage{eufrak}
\usepackage[all]{xy}
\newtheorem{theorem}{Theorem}
\newtheorem{corollary}{Corollary}
\newtheorem{definition}{Definition}
\newtheorem{proposition}{Proposition}

\newtheorem{example}{Example}

\newcommand{\im}{\mathrm{im\,}}

\begin{document}
%\selectlanguage{english}

\title{A SPECTRAL SEQUENCE ASSOCIATED WITH A SYMPLECTIC MANIFOLD}
%  §¡¨¥­¨¥ ­  áâப¨ ®áãé¥á⢫ï¥âáï ª®¬ ­¤®© \\

\author{\firstname{C.}~\surname{Di Pietro}}
\email{cdipietr@unisa.it}
\affiliation{%
Dipartimento di Matematica e Informatica, Universit\`{a} degli
Studi di Salerno, Via Ponte don Melillo, 84084 Fisciano (SA),
Italy}
\author{\firstname{A.~M.}~\surname{Vinogradov}}
% ‡¤¥áì à §¡¨¥­¨¥ ­  áâப¨ ®áãé¥á⢫ï¥âáï  ¢â®¬ â¨ç¥áª¨ ¨«¨ ª®¬ ­¤®© \\
\email{vinograd@unisa.it} \affiliation{ Dipartimento di Matematica
e Informatica, Universit\`{a} degli Studi di Salerno, Via Ponte
don Melillo, 84084 Fisciano (SA), Italy }\affiliation{INFN, Gruppo
Collegato di Salerno, Italy}

\begin{abstract}
With a symplectic manifold a spectral sequence converging to its
de Rham cohomology is associated. A method of computation of its
terms is presented together with some stabilization results. As an
application a characterization of symplectic harmonic manifolds is
given and a relationship with the $\cal{C}$--spectral sequence is
indicated.
\end{abstract}
\maketitle Let $(M,\Omega)$ be a $2n$--dimensional symplectic
manifold and $\Lambda(M)$ be the algebra of differential forms on
$M$. Consider the ideal $\Lambda_{\mathcal{L}}(M)$ of
$\Lambda(M)$, composed of all differential forms that vanish when
restricted to any Lagrangian submanifold of $M$. This ideal is
differentially closed and its powers constitute the
\emph{symplectic filtration} in the de Rham complex of {M}. The
corresponding spectral sequence $ \{E_r^{p,q}, d_r^{p,q} \}$ is
called the \emph{symplectic spectral sequence} associated with
$(M,\Omega)$.

A motivation for this construction comes from the theory of
$\mathcal{C}$--spectral sequences (see \cite{Vin}). Moreover, if
$M= T^*N$, then the symplectic spectral sequence is  nothing but
the ``classical part'' of the $\mathcal{C}$--spectral sequence
associated with the differential equation $d \rho=0$, $\rho \in
\Lambda^1(N)$.

\section{Notations and preliminaries}

In this section the notation is fixed and all necessary facts
concerning symplectic  manifolds (see \cite{Bry,Lyc,Wei} for
further details) are collected.

Throughout the paper $(M, \Omega)$ stands for a $2n$--dimensional
symplectic manifold , $\Lambda = \sum_k \Lambda^k$ for the algebra
of differential forms on $M$, $ H(M)= \sum_k H^k(M)$ for the de
Rham cohomology of $M$ and $D= \sum_k D_k$ for the algebra of
multivectors on $M$.

The isomorphism $\Gamma_1 \colon V \in D_1 \mapsto V \Big\lrcorner\,
\Omega \in \Lambda^1$ of $C^\infty(M)$--modules extends uniquely
to a $C^{\infty}(M)$--algebra isomorphism $\Gamma \colon D \to \Lambda$.
$P= \Gamma^{-1}(\Omega)$ is called the \emph{corresponding to
${\Omega}$ Poisson bivector}.
 $C^{\infty}(M)$--linear operators
$$\top \colon \Lambda^k \to \Lambda^{k+2} \, , \quad  \top \omega = \omega \wedge \Omega, $$
$$\bot \colon \Lambda^k \to \Lambda^{k-2} \, , \quad \bot \omega = P \Big\lrcorner \, \omega, $$
acting on $\Lambda$ are basic for our purposes. Put $\top \Lambda
= \im \top$ and $\Lambda_{\epsilon} = \ker \bot$. Elements of
$\Lambda_\epsilon$ are called \emph{effective forms}.

Another very useful fact is the \emph{Hodge--Lepage} expansion
(see, for instance, \cite{Lyc}) :

\begin{proposition}
Any $\omega\in \Lambda^k, \;k\geq 0,$ admits a unique expansion of
the form
$$\omega = \omega_0 + \top \omega_1 + \top^2 \omega_2 \cdots +
\top^n \omega_n$$ with $\omega_i \in \Lambda_\epsilon^{k-2i}$. In
particular, $\Lambda$ splits into the direct sum $\Lambda_\epsilon
\oplus \top \Lambda$.
\end{proposition}

There exists a symplectic analogue of the Hodge star--operator
$$*\colon\Lambda^k\rightarrow \Lambda^{2n-k}, \;k\geq 0 ,$$
uniquely characterized by the following property: $$\eta \wedge *
\, \omega = (\bot^k (\eta \wedge \omega)) \Omega^n, \quad\eta,
\omega\in\Lambda^k.$$ The map $\delta= \{ \delta_k \}_k$,
$\delta_k = (-1)^{k+1} * d *\colon \Lambda^k \to \Lambda^{k-1}$, is a
$(-1)$--degree differential in $\Lambda$ and $[\bot,d]= \delta$
(see \cite{Bry}).

The introduced operators are subject to the following graded
commutation relations:
\begin{equation} \label{dT}
[d, \delta]=0 \quad ; \quad [\top,d] =0 \quad ;  \quad [\top,\delta]=d \quad ; \quad [\bot, \delta]=0,
\end{equation}
the first of which shows that $(\Lambda, d, \delta)$ is a bicomplex (see \cite{Bry}).

\begin{proposition} \label{deltaef}
The  Hodge--Lepage expansion of $d\omega, \;\omega \in
\Lambda_\epsilon$, is of the form $$d \omega= (d \omega)_0 + \top
(d \omega)_1, \quad (d\omega)_0, \, (d\omega)_1 \in
\Lambda_\epsilon.$$
\end{proposition}
\begin{corollary} \label{Def}
$\Lambda_\epsilon$ is  $\delta$--closed.
\end{corollary}

\section{The term $E_0$}

An explicit description of the symplectic filtration is based on
the following fact of linear algebra.
\begin{proposition} \label{LT}
$\Lambda_{\mathcal{L}} = \top \Lambda$.
\end{proposition}
It shows that the operator $\top^k$ respects the symplectic
filtration by shifting it by $2k$. Moreover, it commute with $d$
and, so, induces an automorphism $\tau=\{\tau_r^k\}, \;\tau_r^k \colon
E_r^{p,q} \to E_r^{p+k,q+k},$ of the symplectic spectral sequence.

\begin{proposition} \label{E0}
${}$
\begin{enumerate}
\item The term $E_0^{p,q}$ is trivial if  $(p,q) \in \Bbb{Z}^2$ lies outside the triangle with
vertexes at $(0,0)$, $(0,n)$ and $(n,n)$.
\item The term $E_0^{0,q}$ is naturally isomorphic to $\Lambda_\epsilon^q$.
\item $\tau_0^p\colon E_0^{0,q-p} \to E_0^{p,q}$ is an isomorphism, if \;$0 \leq p \leq q \leq n$.
\end{enumerate}
\end{proposition}

The last assertion of Propostion \ref{E0} is inductively
generalized to terms $E_r$, $r>0$.

\begin{proposition} \label{tau}
Let $r>0$. Then $\tau_r^k\colon E_r^{p,q} \to E_r^{p+k,q+k}$ is an
isomorphism, if {} $\sum_{i=1}^r (i-1) \leq p \leq q \leq q+k \leq
n - (2+\sum_{i=1}^r(i-2))$. In particular, $\tau_1^p\colon E_1^{0,q-p}
\to E_1^{p,q}$ is an isomorphism, if \;$0 \leq p \leq q < n$.
\end{proposition}

\section{The term $E_1$}

The exact sequence constructed in this section gives a useful
description of the first term of the symplectic spectral sequence.
It is composed of the following two families of
$\Bbb{R}$--homomorphisms;
\begin{align*}
\phi_{p,q}\colon &  H^{q-p}(M)   \to  E_1^{p,q}, \qquad 0 \leq p \leq q \leq n\\
 & [\omega]_{\im d}   \mapsto [\top^p \omega ]_{\im d_0} \\
\psi_{p,q}\colon& E_1^{p,q}  \to H^{q-p-1}(M), \qquad 0 \leq p \leq q < n\\
 & [\top^p \rho]_{\im d_0} \mapsto [\eta]_{\im d}
\end{align*}
with $\top \eta = d \rho$.

Put $A^k = \{\omega \ | \ \, \omega \in \Lambda^k ; \, d \omega
\in \Lambda_\epsilon^{k+1} \}$, $C^{k} = A^k / \im d \cap
\Lambda^k$ and consider $\Bbb{R}$--homomorphisms
\begin{align*}
\phi_{p,n} \colon  C^{n-p}  \to  E_1^{p,n}, \quad & [\omega]_{\im d} \mapsto [\top^p \omega]_{\im d_0} \\
\psi_{p,n} \colon  E_1^{p,n}  \to C^{n-(p+1)}, \quad & [\top^p \rho]_{\im d_0} \mapsto [\eta]_{\im d}.
\end{align*}
The above defined homomorphisms together with the multiplication
by the cohomology class of $\Omega$ homomorphism
$\tau \colon H^*(M)\rightarrow H^{*+2}(M)$ form the following sequences
\begin{equation}\label{LES}
\dots\rightarrow H^{q-p}(M) \stackrel{\phi_{p,q}}{\longrightarrow}
 E_1^{p,q}\stackrel{\psi_{p,q}}{\longrightarrow}
H^{q-(p+1)}(M)\stackrel{\tau}{\longrightarrow}H^{q+1-p}(M)\stackrel
{\phi_{p,q+1}} {\longrightarrow} E_1^{p,q+1}\rightarrow\dots
\end{equation}
whose left and right ends are
\begin{equation*}
0\rightarrow H^0(M)\stackrel{\phi_{p,p}}{\longrightarrow}
 E_1^{p,p}\stackrel{\psi_{p,p}}{\longrightarrow}0\stackrel{\tau}{\longrightarrow}
H^1(M)\stackrel{\phi_{p,p+1}}{\longrightarrow}
 E_1^{p,p+1}\rightarrow\dots
\end{equation*}
and
\begin{equation*}
\dots\rightarrow H^{n-(p+2)}(M)\stackrel{\tau}{\longrightarrow}
C^{n-p}\stackrel{\phi_{p,n}}{\longrightarrow}E_1^{p,n}
\stackrel{\psi_{p,n}}{\longrightarrow}C^{n-(p+1)}
\stackrel{\tau^{p+1}}{\longrightarrow}H^{n+p+1}(M)\rightarrow 0\;,
\end{equation*}
respectively. Notice that sequence (\ref{LES}) involves the terms
of the $p$--th column of $E_1$.

\begin{theorem} \label{TLES}
 The sequence (\ref{LES}) is exact.
\end{theorem}

\section{Stabilization theorems}

Theorem \ref{TLES} is key in studying stability of the symplectic
spectral sequence. For instance, if $\Omega$ is exact, then the
homomorphism $\tau\colon H^k(M) \to H^{k+2}(M)$ in (\ref{LES}) is
trivial. This fact and Theorem \ref{TLES} lead to the following
result.

\begin{theorem} \label{Oex}
If $\Omega$ is exact, then the symplectic spectral sequence
stabilizes at the term $E_2$. Moreover, if a term $E_2^{p,q}$ is
different from zero, then either $p=0$, or $q=n$ and $
E_2^{0,q}\cong H^{q}(M)$, \;if $q \leq n$, and $ E_2^{p,n}\cong
H^{n+p}(M)$, \;if $p \geq 0$.
\end{theorem}

\begin{corollary}
If $\Omega$ is the standard symplectic form on $M=T^*N$, $\mathrm{dim} \, N=n$,
then the corresponding symplectic spectral sequence stabilizes at the second term
and $E_2^{0,q} \cong H^q(N)$, $E_2^{p,q}=0$, for $p>0$.
\end{corollary}

Theorem \ref{Oex} is generalized as follows.

\begin{theorem} \label{EST}
Let $t > 1$ be the minimal integer such that $\Omega^t = d \rho$.
Then the symplectic spectral sequence stabilizes at the term
$E_{t-r+1}$ where $r$ is an integer such that $0 \leq r \leq
\mathrm{max} \{ 0, 2t-(n+1) \}$.
\end{theorem}

The estimate in the previous theorem can not be, generally,
improved as the following example shows .

\begin{example}
The symplectic spectral sequence associated with $(C, \Omega)
\times (\Bbb{R}^{2h}, \Omega_{\Bbb{R}^{2h}})$, where $(C, \Omega)$
is a closed symplectic $2m$--fold and $(\Bbb{R}^{2h},
\Omega_{\Bbb{R}^{2h}})$ is the standard symplectic manifold, is
not stable in terms $E_{h}$, if $h \leq m+1=t$, and becomes stable
by starting from the term $E_{m+2}$.
\end{example}

If cohomology classes $[\Omega^t], \;t \leq n$, are all
nontrivial, then $M$ is closed. Closed symplectic manifolds are
characterized by the fact that $E_\infty^{p,p}= E_2^{p,p} \cong
\mathbb{R}$, if $0 \leq p \leq n$. Moreover, it holds

\begin{theorem}
The symplectic spectral sequence for a closed symplectic manifold
stabilizes at the term $E_2$.
\end{theorem}

\subsection{The Brylinski conjecture}

\begin{definition}
A form $\eta$ is called \emph{symplectically harmonic} iff $\eta
\in \ker d \cap \ker \delta$.
\end{definition}

In \cite{Bry} Brylinski conjectured that each cohomology class of
a closed symplectic manifold contains at least one symplectically
harmonic form. It was, however, disproved by counterexamples found
by Mathieu (see \cite{Mat}). So, the problem of characterization
of simplectic manifolds for which the Brylinski conjecture holds
arises. Call such symplectic manifolds \emph{harmonic}. In view of
(\ref{dT}) and Proposition \ref{deltaef} we have the following
useful technical characterization of harmonic symplectic
manifolds.

\begin{proposition}
A symplectic manifolds is harmonic iff any cohomology class is
represented by a form $\eta$ whose Hodge--Lepage expansion terms
are all closed.
\end{proposition}

This Proposition helps to prove

\begin{theorem}
A closed symplectic manifold is harmonic  iff $\tau_2^{n-q} \colon E_2^{0,q} \to E_2^{n-q,n}$
is an isomorphism for all $0
\leq q \leq n$.
\end{theorem}

\begin{corollary}
For a closed harmonic manifold mappings $\tau_2^{n-p-q}\colon E_2^{p,q} \to E_2^{n-q, n-p}$
are isomorphisms for $p+q \leq n$ and $0 \leq p \leq q$. In particular, $E_2$ is symmetric
with respect to the line $p+q=n$.
\end{corollary}

\section{Associated diffiety}

Let  $J^k(M,n)$ be the manifold of $k$--jets of $n$--dimensional
submanifolds of a manifold $M$. The \emph{$k$--th prolongation
$M_{(k)}$} of a symplectic manifold $(M,\Omega)$ is a submanifold
in $J^k(M,n)$ composed of $k$--th jets of Lagrangian submanifolds
of $M$. By restricting the Cartan distribution on $J^\infty(M,n)$
to $M_{(\infty)}$ we obtain a diffiety which locally coincides
with the infinite prolongation of the equation $d\rho=0, \rho\in
\Lambda^1(N), \;\dim N=n$. A natural projection
$M_{(\infty)}\rightarrow M$ induces a morphism of the symplectic
spectral sequence of $M$ to the $\cal{C}$--spectral sequence of
the diffiety $M_{(\infty)}$. This allows to find out some useful
interpretations for various terms and differentials of the
symplectic spectral sequence. For instance, consider the following
functional defined on compact Lagrangian submanifolds associated
with an element $\theta\in E_1^{0,n}$:
\begin{equation*}
L \quad\mapsto \int_L\omega|_L, \qquad \theta=[\omega],
\qquad\omega\in \Lambda^n \mod \top\Lambda.
\end{equation*}
Then $d^{0,n}_1(\theta)=0$ is the Euler-Lagrange equation for
extremals of this functional.

\section{Generalisations}

There are numerous analogues of the symplectic spectral sequence
due to the  fact that the underling construction is of a rather
general nature. Below we list some of its "neighbors".
\begin{itemize}
\item Let $N$ be a submanifold of a symplectic manifold $M$. By
restricting the symplectic filtration in $\Lambda(M)$ to
$\Lambda(N)$ we get a spectral sequence converging to the de Rham
cohomology of $N$. In particular, this way a spectral sequence is
associated with a Hamilton--Jacobi equation $\mathcal E$. This
spectral sequence is the "classical part" of the $\mathcal
C$--spectral sequence associated with $\mathcal E$.
\item Let $M$ be a contact manifold and
$\mathcal{I}\subset\Lambda(M)$ the ideal composed of differential
forms that vanish on all Legendre submanifolds of $M$. Powers of
$\mathcal{I}$ form the \emph{contact filtration} in $\Lambda(M)$.
This way one gets the \emph{contact spectral sequence} associated
with a contact manifold.
\item By restricting the contact filtration to a submanifold $N$ of a
contact manifold $M$ one gets a spectral sequence which is the
classical part of the $\mathcal C$--spectral sequence associated
with $N$ interpreted as an (overdetermined) system of first order
scalar differential equations (see \cite{Lyc,Vin}).
\item Let $(M,P)$ be a Poisson manifold, $P$ being the Poisson bivector.
The \emph{Poisson filtration} in the algebra $D(M)$ of
multivectors fields on $M$ is that formed by powers of the
principal ideal generated by $P$. The Poisson differential $d_P,
\;d_P(Q)=[[P,Q]], \;Q\in D(M),$ respects this filtration and the
corresponding to it spectral sequence is called \emph{Poisson}. If
$P$ is nondegerate the Poisson spectral sequence is, in a sense,
dual to the corresponding symplectic one.
\end{itemize}

Constructions and results of this note naturally generalize to all
these spectral sequences. Details will be given in a separate
publication.

\end{document}